\documentclass[12pt]{article}

\newcommand{\qed}{{\hfill$\square$}}

\usepackage{amsmath, epsfig, cite, lineno}
\usepackage{amssymb}
\usepackage{amsfonts}
\usepackage{latexsym}

\newtheorem{thm}{Theorem}[section]

\newtheorem{cor}[thm]{Corollary}

\newtheorem{lem}[thm]{Lemma}
\newtheorem{conj}[thm]{Conjecture}

\def\pf{\noindent {\it Proof.} }

\def\fl#1{\left\lfloor#1\right\rfloor}

\numberwithin{equation}{section}

\begin{document}


\begin{center}
{\Large\bf Proof of two divisibility properties of binomial\\[5pt] coefficients conjectured by Z.-W. Sun}
\end{center}

\vskip 2mm \centerline{Victor J. W. Guo}
\begin{center}
{\footnotesize Department of Mathematics, Shanghai Key Laboratory of PMMP,
East China Normal University,\\ 500 Dongchuan Rd., Shanghai 200241,
 People's Republic of China\\
{\tt jwguo@math.ecnu.edu.cn,\quad
http://math.ecnu.edu.cn/\textasciitilde{jwguo} }}
\end{center}

\vskip 0.7cm
{\small \noindent{\bf Abstract.} For all positive integers $n$, we prove the following divisibility properties:
\begin{align*}
(2n+3){2n\choose n}  \left|3{6n\choose 3n}{3n\choose n},\right. \quad\text{and}\quad
(10n+3){3n\choose n} \left|21{15n\choose 5n}{5n\choose n}.\right.
\end{align*}
This confirms two recent conjectures of Z.-W. Sun. Some similar divisibility properties are given.
Moreover, we show that, for all positive integers $m$ and $n$, the product
$am{am+bm-1\choose am}{an+bn\choose an}$ is divisible by $m+n$.
In fact, the latter result can be generalized to the $q$-binomial coefficients and $q$-integers case,
which generalizes the positivity of $q$-Catalan numbers. We also propose several related conjectures.
}



\section{Introduction}
In \cite{Sun3,Sun2}, Z.-W. Sun proved some divisibility properties of binomial coefficients, such as
\begin{align*}
2(2n+1){2n\choose n} &\left|{6n\choose 3n}{3n\choose n}\right.,\\[5pt]    
(10n+1){3n\choose n} &\left|{15n\choose 5n}{5n-1\choose n-1}.\right.
\end{align*}
Some similar divisibility results were later obtained by Guo \cite{Guo} and Guo and Krattenthaler \cite{GK}.
It should be mentioned that Bober \cite{Bober} has completely described
when ratios of factorial products of the form
$$
\frac{(a_1 n)!\cdots (a_{k}n)!}{(b_1 n)!\cdots (b_{k+1}n)!}
$$
with $a_1+\cdots+a_k=b_1+\cdots+b_{k+1}$ are always integers.

Let
\begin{align*}
S_n=\frac{{6n\choose 3n}{3n\choose n}}{2(2n+1){2n\choose n}},\quad\text{and}\quad
t_n=\frac{{15n\choose 5n}{5n-1\choose n-1}}{(10n+1){3n\choose n}}.
\end{align*}
In this paper we first prove the following two results conjectured by Z.-W. Sun \cite{Sun3,Sun2}.
\begin{thm}\label{conj:sun3}
{\rm(see \cite[Conjecture 3(i)]{Sun3})}Let $n$ be a positive
integer. Then
\begin{align*}
3 S_n \equiv 0 \pmod{2n+3}.
\end{align*}
\end{thm}
\begin{thm}\label{conj:sun2}
{\rm\cite[Conjecture 1.3]{Sun2}}Let $n$ be a positive integer. Then
\begin{align*}
21 t_n \equiv 0 \pmod{10n+3}.
\end{align*}
\end{thm}
We shall also give more congruences for $S_n$ and $t_n$ as follows.
\begin{thm}\label{thm:new1}
Let $n$ be a positive integer. Then
\begin{align}
105 S_n &\equiv 0 \pmod{2n+5}, \label{eq:number1} \\
315 S_n &\equiv 0 \pmod{2n+7}, \label{eq:number2}\\
6435 S_n &\equiv 0 \pmod{2n+9}, \\
3003 t_n & \equiv 0 \pmod{2n+1}, \\
88179 t_n & \equiv 0 \pmod{10n+7}, \\
43263 t_n & \equiv 0 \pmod{10n+9}. \label{eq:number6}
\end{align}
\end{thm}
Let $\mathbb{Z}$ denote the set of integers. Another result in this paper is the following.
\begin{thm}\label{thm:new2}
Let $a,b,m,n$ be positive integers. Then
\begin{align}
\frac{abm}{(a+b)(m+n)}{am+bm\choose am}{an+bn\choose an}
=\frac{am}{m+n}{am+bm-1\choose am}{an+bn\choose an}
\in\mathbb{Z}. \label{eq:abmn}
\end{align}
\end{thm}
Letting $a=b=1$ in \eqref{eq:abmn}, we get the following result.
\begin{cor}
Let $m,n$ be positive integers. Then
\begin{align}
\frac{m}{2(m+n)}{2m\choose m}{2n\choose n}\in\mathbb{Z}.
\end{align}
In particular,
\begin{align*}
\frac{6}{n+2}{2n\choose n},\quad \frac{30}{n+3}{2n\choose n},\quad
\frac{140}{n+4}{2n\choose n},\quad \frac{630}{n+5}{2n\choose n}\in\mathbb{Z}.
\end{align*}
\end{cor}

In the next section, we give some lemmas. The proofs of Theorems \ref{conj:sun3}--\ref{thm:new1}
will be given in Sections 3--5 respectively. A proof of the $q$-analogue of Theorem \ref{thm:new2}
will be given in Section 6. We close our paper with some further remarks and open problems in Section 7.

\section{Some lemmas}
For the $p$-adic order of $n!$, there is a known formula
\begin{align}
{\rm ord}_p n!=\sum_{i=1}^\infty\left\lfloor\frac{n}{p^i}\right\rfloor,
\label{eq:ord}
\end{align}
where $\lfloor x\rfloor$ denotes the greatest integer not exceeding
$x$. In this section, we give some results on the floor
function $\lfloor x\rfloor$.
\begin{lem}\label{lem:one}
For any real number $x$, we have
\begin{align}
\left\lfloor 6x\right\rfloor+\left\lfloor x\right\rfloor
&\geqslant
\left\lfloor 3x\right\rfloor+2\left\lfloor 2x\right\rfloor, \label{eq:000}\\
\left\lfloor 15x\right\rfloor+\left\lfloor 2x\right\rfloor
&\geqslant
\left\lfloor 10x\right\rfloor+\left\lfloor 4x\right\rfloor
+\left\lfloor 3x\right\rfloor.
\end{align}
\end{lem}
\pf See \cite[Theorem 1.1]{Bober} and one of the 52 sporadic step
functions given in \cite[Table 2, line\# 32]{Bober}.
\qed
\begin{lem}\label{lem:two}
Let $m$ and $n$ be two positive integers such that $m|2n+3$ and $m\geqslant 5$. Then
\begin{align}
\left\lfloor\frac{6n}{m}\right\rfloor+\left\lfloor\frac{n}{m}\right\rfloor
=\left\lfloor\frac{3n}{m}\right\rfloor+2\left\lfloor\frac{2n}{m}\right\rfloor+1.  \label{eq:ineq0}
\end{align}
\end{lem}
\pf Let $\{x\}=x-\lfloor x\rfloor$ be the fractional part of $x$.
Then \eqref{eq:ineq0} is equivalent to
\begin{align}
\left\{\frac{6n}{m}\right\}+\left\{\frac{n}{m}\right\}
=\left\{\frac{3n}{m}\right\}+2\left\{\frac{2n}{m}\right\}-1.  \label{eq:ineq02}
\end{align}
Now suppose that $m|2n+3$ and $m\geqslant 5$. We have
\begin{align*}
\left\{\frac{2n}{m}\right\}=\frac{m-3}{m}>\frac{1}{3},\quad{\rm and}\quad
\left\lfloor\frac{2n}{m}\right\rfloor=\frac{2n+3}{m}-1\equiv
0,2,4,6,8\pmod{10}.
\end{align*}
It follows that
\begin{align*}
\left\{\frac{6n}{m}\right\}
&=\begin{cases}
\displaystyle\frac{2m-9}{m}, &\text{if $m=5,7$}, \\[15pt]
\displaystyle\frac{m-9}{m}, &\text{if $m\geqslant 9$,}
\end{cases} \\[5pt]
\left\{\frac{n}{m}\right\}
&=\frac{m-3}{2m}, \\[5pt]
\left\{\frac{3n}{m}\right\}
&=\begin{cases}
\displaystyle\frac{3m-9}{2m}, &\text{if $m=5,7$}, \\[15pt]
\displaystyle\frac{m-9}{2m}, &\text{if $m\geqslant 9$.}
\end{cases}
\end{align*}
Therefore, the identity \eqref{eq:ineq02} is true for any positive integer $m\geqslant 5$. \qed
\begin{lem}\label{lem:three}
Let $m$ and $n$ be two positive integers such that $m|10n+3$ and
$m\geqslant 9$. Then
\begin{align}
\left\lfloor\frac{15n}{m}\right\rfloor+\left\lfloor\frac{2n}{m}\right\rfloor
=\left\lfloor\frac{10n}{m}\right\rfloor+\left\lfloor\frac{4n}{m}\right\rfloor
+\left\lfloor\frac{3n}{m}\right\rfloor+1.  \label{eq:ineq1}
\end{align}
\end{lem}
\pf It is easy to see that \eqref{eq:ineq1} is equivalent to
\begin{align}
\left\{\frac{15n}{m}\right\}+\left\{\frac{2n}{m}\right\}
=\left\{\frac{10n}{m}\right\}+\left\{\frac{4n}{m}\right\}
+\left\{\frac{3n}{m}\right\}-1.  \label{eq:ineq2}
\end{align}
Now suppose that $m|10n+3$ and $m\geqslant 9$. We have
\begin{align*}
\left\{\frac{10n}{m}\right\}=\frac{m-3}{m}\geqslant \frac{2}{3},\quad{\rm and}\quad
A:=\left\lfloor\frac{10n}{m}\right\rfloor=\frac{10n+3}{m}-1\equiv
0,2,6,8\pmod{10}.
\end{align*}
It is easy to check that
\begin{align*}
\left\{\frac{15n}{m}\right\}&=\frac{m-9}{2m},\\[5pt]
\left(\left\{\frac{2n}{m}\right\},\left\{\frac{4n}{m}\right\},
\left\{\frac{3n}{m}\right\} \right)
& =\begin{cases}
\displaystyle\left(\frac{2m-6}{10m},\frac{4m-12}{10m},\frac{3m-9}{10m}\right),
&\text{if $A\equiv 0\pmod{10},$} \\[15pt]
\displaystyle\left(\frac{6m-6}{10m},\frac{2m-12}{10m},\frac{9m-9}{10m}\right),
&\text{if $A\equiv 2\pmod{10},$} \\[15pt]
\displaystyle\left(\frac{4m-6}{10m},\frac{8m-12}{10m},\frac{m-9}{10m}\right),
&\text{if $A\equiv 6\pmod{10},$} \\[15pt]
\displaystyle\left(\frac{8m-6}{10m},\frac{6m-12}{10m},\frac{7m-9}{10m}\right),
&\text{if $A\equiv 8\pmod{10},$}
\end{cases}
\end{align*}
and so the identity \eqref{eq:ineq2} holds. \qed

\section{Proof of Theorem \ref{conj:sun3}}
Let $\gcd(a,b)$ denote the greatest common divisor of two integers
$a$ and $b$. For any positive integer $n$, since $\gcd(2n+3,4n+2)=1$, to prove  Theorem \ref{conj:sun3}, it is enough
to show that
\begin{align}
(2n+3)\left|\frac{3{6n\choose 3n}{3n\choose n}}{{2n\choose n} }.\right.  \label{eq:div01}
\end{align}
By \eqref{eq:ord}, for any odd prime $p$, the $p$-adic order of
$$
\frac{{6n\choose 3n}{3n\choose n}}{(2n+3){2n\choose n}}=\frac{(2n+2)!(6n)!(n)!}{(2n+3)!(3n)!(2n)!^2}
$$
is given by
\begin{align}
\sum_{i=1}^\infty\left(\left\lfloor\frac{2n+2}{p^i}\right\rfloor
+\left\lfloor\frac{6n}{p^i}\right\rfloor+\left\lfloor\frac{n}{p^i}\right\rfloor
-\left\lfloor\frac{2n+3}{p^i}\right\rfloor-\left\lfloor\frac{3n}{p^i}\right\rfloor
-2\left\lfloor\frac{2n}{p^i}\right\rfloor\right). \label{eq:sum0}
\end{align}
Note that
$$
\left\lfloor\frac{2n+2}{p^i}\right\rfloor
-\left\lfloor\frac{2n+3}{p^i}\right\rfloor
=\begin{cases}
-1,&\text{if $p^i|2n+3,$}\\[5pt]
0,&\text{otherwise.}
\end{cases}
$$
By Lemmas \ref{lem:one} and \ref{lem:two}, for $p\geqslant 5$, the summation \eqref{eq:sum0} is clearly greater than or
equal to $0$. For $p=3$, we have $\eqref{eq:sum0}\geqslant -1$
because if the number $i$ satisfies $3^i|2n+3$ and $3^i<5$ then we must have $i=1$.
This proves that
\begin{align*}
\frac{3{6n\choose 3n}{3n\choose n}}{(2n+3){2n\choose n}}
\end{align*}
is always an integer. Hence \eqref{eq:div01} holds.

\medskip\noindent{\it Remark.} Z.-W. Sun \cite[Conjecture 3(i)]{Sun3} also conjectured that
$S_n$ is odd if and only if $n$ is a power of $2$. After reading a previous version of this paper,
Quan-Hui Yang told me that it is easy to show that
${\rm ord}_2((6n)!n!/(3n)!(2n)!^2)$ equals the number of $1$'s in the binary expansion of $n$ by
noticing that
$$
{\rm ord}_2 (6n)!=3n+{\rm ord}_2 (3n)!,\quad {\rm ord}_2 (2n)!=n+{\rm ord}_2 n!,
$$
and using Legendre's theorem. T. Amdeberhan and V.H. Moll also pointed out this.

\section{Proof of Theorem \ref{conj:sun2}}
For any positive integer $n$, since $\gcd(10n+1,10n+3)=1$, to prove Theorem \ref{conj:sun2},
it is enough to show that
\begin{align}
(10n+3)\left|\frac{21{15n\choose 5n}{5n-1\choose n-1}}{{3n\choose n} }.\right.  \label{eq:div1}
\end{align}
Furthermore, since $\gcd(10n+1,5)=1$ and ${5n\choose n}=5{5n-1\choose n-1}$,
one sees that  \eqref{eq:div1} is equivalent to
\begin{align}
(10n+3)\left|\frac{21{15n\choose 5n}{5n\choose n}}{{3n\choose n} }.\right.  \label{eq:div2}
\end{align}
By \eqref{eq:ord}, for any odd prime $p$, the $p$-adic order of
$$
\frac{{15n\choose 5n}{5n\choose n}}{(10n+3){3n\choose n}}=\frac{(10n+2)!(15n)!(2n)!}{(10n+3)!(10n)!(4n)!(3n)!}
$$
is given by
\begin{align}
\sum_{i=1}^\infty\left(\left\lfloor\frac{10n+2}{p^i}\right\rfloor
+\left\lfloor\frac{15n}{p^i}\right\rfloor+\left\lfloor\frac{2n}{p^i}\right\rfloor
-\left\lfloor\frac{10n+3}{p^i}\right\rfloor-\left\lfloor\frac{10n}{p^i}\right\rfloor
-\left\lfloor\frac{4n}{p^i}\right\rfloor-\left\lfloor\frac{3n}{p^i}\right\rfloor\right). \label{eq:sum}
\end{align}
Note that
$$
\left\lfloor\frac{10n+2}{p^i}\right\rfloor
-\left\lfloor\frac{10n+3}{p^i}\right\rfloor
=\begin{cases}
-1,&\text{if $p^i|10n+3,$}\\[5pt]
0,&\text{otherwise.}
\end{cases}
$$
By Lemmas \ref{lem:one} and \ref{lem:three}, for $p\geqslant 11$, or $p=5$, the summation \eqref{eq:sum} is clearly greater than or
equal to $0$. For $p=3,7$, we have $\eqref{eq:sum}\geqslant -1$ because there is at most one index $i$
satisfying $p^i|10n+3$ and $p^i<9$ in this case. This proves that
\begin{align*}
\frac{21{15n\choose 5n}{5n\choose n}}{(10n+3){3n\choose n}}
\end{align*}
is always an integer. Namely, \eqref{eq:div2} is true.

\section{Proof of Theorem \ref{thm:new1}}
\begin{lem}\label{lem:4}
Let $m$ and $n$ be two positive integers. Then
\begin{align}
\left\lfloor\frac{6n}{m}\right\rfloor+\left\lfloor\frac{n}{m}\right\rfloor
=\left\lfloor\frac{3n}{m}\right\rfloor+2\left\lfloor\frac{2n}{m}\right\rfloor+1,  \label{eq:lem4}
\end{align}
if $m|2n+5$ and $m\geqslant 9$, or $m|2n+7$ and $m\geqslant 11$,  or $m|2n+9$ and $m\geqslant 15$.
\end{lem}
\pf The proof is similar to that of Lemma \ref{lem:two}.
We only consider the case when $m|2n+5$ and $m\geqslant 9$. In this case, we have
\begin{align*}
\left\{\frac{2n}{m}\right\}=\frac{m-5}{m}>\frac{1}{3},\quad{\rm and}\quad
\left\lfloor\frac{2n}{m}\right\rfloor=\frac{2n+5}{m}-1\equiv
0,2,4,6,8\pmod{10}.
\end{align*}
It follows that
\begin{align*}
\left\{\frac{6n}{m}\right\}
&=\begin{cases}
\displaystyle\frac{2m-15}{m}, &\text{if $m=9,11,13$}, \\[15pt]
\displaystyle\frac{m-15}{m}, &\text{if $m\geqslant 15$,}
\end{cases} \\[5pt]
\left\{\frac{n}{m}\right\}
&=\frac{m-5}{2m}, \\[5pt]
\left\{\frac{3n}{m}\right\}
&=\begin{cases}
\displaystyle\frac{3m-15}{2m}, &\text{if $m=9,11,13$}, \\[15pt]
\displaystyle\frac{m-15}{2m}, &\text{if $m\geqslant 15$.}
\end{cases}
\end{align*}
Therefore,
\begin{align*}
\left\{\frac{6n}{m}\right\}+\left\{\frac{n}{m}\right\}
=\left\{\frac{3n}{m}\right\}+2\left\{\frac{2n}{m}\right\}-1,
\end{align*}
and so \eqref{eq:lem4} holds.
\qed

\begin{lem}\label{lem:5}
Let $m$ and $n$ be two positive integers. Then
\begin{align}
\left\lfloor\frac{15n}{m}\right\rfloor+\left\lfloor\frac{2n}{m}\right\rfloor
=\left\lfloor\frac{10n}{m}\right\rfloor+\left\lfloor\frac{4n}{m}\right\rfloor
+\left\lfloor\frac{3n}{m}\right\rfloor+1, \label{eq:lem5}
\end{align}
if $m|2n+1$ and $m\geqslant 15$, or $m|10n+7$ and $m\geqslant 21$,  or $m|10n+9$ and $m\geqslant 27$.
\end{lem}
\pf The proof is similar to that of Lemma \ref{lem:three}.
We only consider the case when $m|10n+9$ and $m\geqslant 27$. In this case, we have
\begin{align*}
\left\{\frac{10n}{m}\right\}=\frac{m-9}{m}\geqslant \frac{2}{3},\quad{\rm and}\quad
A:=\left\lfloor\frac{10n}{m}\right\rfloor=\frac{10n+9}{m}-1\equiv
0,2,6,8\pmod{10}.
\end{align*}
It follows that
\begin{align*}
\left\{\frac{15n}{m}\right\}&=\frac{m-27}{2m},\\[5pt]
\left(\left\{\frac{2n}{m}\right\},\left\{\frac{4n}{m}\right\},
\left\{\frac{3n}{m}\right\} \right)
& =\begin{cases}
\displaystyle\left(\frac{2m-18}{10m},\frac{4m-36}{10m},\frac{3m-27}{10m}\right),
&\text{if $A\equiv 0\pmod{10},$} \\[15pt]
\displaystyle\left(\frac{6m-18}{10m},\frac{2m-36}{10m},\frac{9m-27}{10m}\right),
&\text{if $A\equiv 2\pmod{10},$} \\[15pt]
\displaystyle\left(\frac{4m-18}{10m},\frac{8m-36}{10m},\frac{m-27}{10m}\right),
&\text{if $A\equiv 6\pmod{10},$} \\[15pt]
\displaystyle\left(\frac{8m-18}{10m},\frac{6m-36}{10m},\frac{7m-27}{10m}\right),
&\text{if $A\equiv 8\pmod{10}.$}
\end{cases}
\end{align*}
Hence,
\begin{align*}
\left\{\frac{15n}{m}\right\}+\left\{\frac{2n}{m}\right\}
=\left\{\frac{10n}{m}\right\}+\left\{\frac{4n}{m}\right\}
+\left\{\frac{3n}{m}\right\}-1,
\end{align*}
and \eqref{eq:lem5} holds.
\qed

\medskip
\noindent{\it Proof of Theorem {\rm\ref{thm:new1}.}} Since the proofs of
the congruences \eqref{eq:number1}--\eqref{eq:number6} are similar in view of Lemmas \ref{lem:4} and \ref{lem:5},
we only give proofs of \eqref{eq:number2} and \eqref{eq:number6}. Noticing that
$\gcd(2n+1,2n+7)=1$ or $3$, to prove \eqref{eq:number2}, it suffices to show that
\begin{align}
(2n+7)\left|\frac{105{6n\choose 3n}{3n\choose n}}{{2n\choose n} }.\right.  \label{eq:div04}
\end{align}
Let
$$
X_n:=\frac{{6n\choose 3n}{3n\choose n}}{(2n+7){2n\choose n}}=\frac{(2n+6)!(6n)!(n)!}{(2n+7)!(3n)!(2n)!^2}.
$$
By \eqref{eq:ord}, for any odd prime $p$, we have
\begin{align*}
{\rm ord}_p X_n=\sum_{i=1}^\infty\left(\left\lfloor\frac{2n+6}{p^i}\right\rfloor
+\left\lfloor\frac{6n}{p^i}\right\rfloor+\left\lfloor\frac{n}{p^i}\right\rfloor
-\left\lfloor\frac{2n+7}{p^i}\right\rfloor-\left\lfloor\frac{3n}{p^i}\right\rfloor
-2\left\lfloor\frac{2n}{p^i}\right\rfloor\right).
\end{align*}
Note that \eqref{eq:lem4} is also true for $m=3$ and $n\equiv 1\pmod 3$, and
$$
\left\lfloor\frac{2n+6}{p^i}\right\rfloor
-\left\lfloor\frac{2n+7}{p^i}\right\rfloor
=\begin{cases}
-1,&\text{if $p^i|2n+7,$}\\[5pt]
0,&\text{otherwise.}
\end{cases}
$$
By Lemmas \ref{lem:one} and \ref{lem:4}, we obtain
\begin{align*}
\begin{cases}
{\rm ord}_p X_n\geqslant 0,&\text{if $p\geqslant 11$},\\[5pt]
{\rm ord}_p X_n\geqslant -1,&\text{if  $p=3,5,7$}.
\end{cases}
\end{align*}
This proves \eqref{eq:div04}.

Similarly, since $\gcd(10n+1,10n+9)=\gcd(10n+9,5)=1$, the congruence \eqref{eq:number6}
is equivalent to
\begin{align}
(10n+9)\left|\frac{43263{15n\choose 5n}{5n\choose n}}{{3n\choose n} }.\right.  \label{eq:div05}
\end{align}
Let
$$
Y_n:=\frac{{15n\choose 5n}{5n\choose n}}{(10n+9){3n\choose n}}
=\frac{(10n+8)!(15n)!(2n)!}{(10n+9)!(10n)!(4n)!(3n)!}
$$
Then ${\rm ord}_p Y_n$ is given by
\begin{align*}
\sum_{i=1}^\infty\left(\left\lfloor\frac{10n+8}{p^i}\right\rfloor
+\left\lfloor\frac{15n}{p^i}\right\rfloor+\left\lfloor\frac{2n}{p^i}\right\rfloor
-\left\lfloor\frac{10n+9}{p^i}\right\rfloor-\left\lfloor\frac{10n}{p^i}\right\rfloor
-\left\lfloor\frac{4n}{p^i}\right\rfloor-\left\lfloor\frac{3n}{p^i}\right\rfloor\right).
\end{align*}
Note that \eqref{eq:lem5} also holds for  $m|10n+7$ and $m=7,13,17$.
Similarly as before, we have
\begin{align*}
\begin{cases}
{\rm ord}_p Y_n\geqslant 0,&\text{if $p=5,7,13,17$, or $p\geqslant 29,$}\\[5pt]
{\rm ord}_p Y_n\geqslant -1,&\text{if  $p=11,19,23,$}\\[5pt]
{\rm ord}_p Y_n\geqslant -2,&\text{if  $p=3.$}
\end{cases}
\end{align*}
Observing that $43263=3^2\cdot 11\cdot 19\cdot 23$, we complete the proof of \eqref{eq:div05}.

\section{A $q$-analogue of Theorem \ref{thm:new2}}
Recall that the {\it $q$-binomial coefficients}
are defined by
\begin{align*}
{n\brack k}_q
=\begin{cases}
\displaystyle\frac{(1-q^n)(1-q^{n-1})\cdots (1-q^{n-k+1})}{(1-q)(1-q^2)\cdots(1-q^k)},
&\text{if $0\leqslant k\leqslant n,$} \\[5pt]
0, &\text{otherwise.}
\end{cases}
\end{align*}
We begin with the announced strengthening of Theorem \ref{thm:new2}. It is easily seen that Theorem \ref{thm:new2}
can be obtained upon letting $q\to 1$ in Corollary \ref{cor:2}.
\begin{thm}\label{thm:final}
Let $a,b,m,n\geqslant 1$. Then
\begin{align}
\frac{1-q^{\gcd(am,m+n)}}{1-q^{m+n}}{am+bm-1\brack am}_q {an+bn\brack an}_q  \label{eq:last}
\end{align}
is a polynomial in $q$ with non-negative integer coefficients.
\end{thm}

\begin{cor}\label{cor:2}
Let $a,b,m,n\geqslant 1$. Then
\begin{align}
\frac{1-q^{am}}{1-q^{m+n}}{am+bm-1\brack am}_q {an+bn\brack an}_q \label{eq:last00}
\end{align}
is a polynomial in $q$ with non-negative integer coefficients.
\end{cor}

It is clear that, when $a=b=m=1$, the numbers \eqref{eq:last00} reduce
to the $q$-Catalan numbers
$$
C_n(q)=\frac{1-q}{1-q^{2n+1}}{2n\brack n}_q.
$$
It is well known that the $q$-Catalan numbers $C_n(q)$ are
polynomials with non-negative integer coefficients
(see \cite{Andrews87,Andrews93,Andrews10,FH}).
There are many different $q$-analogues
of the Catalan numbers (see F\"{u}rlinger and Hofbauer \cite{FH}).
For the so-called $q,t$-Catalan numbers, see \cite{GH,Haiman,Haglund}.

Recall that a  polynomial
$P(q)=
\sum _{i=0} ^{d}p_iq^i$ in $q$ of degree $d$ is called {\it reciprocal\/} if
$p_i=p_{d-i}$ for all $i$, and that it
is called {\it unimodal\/} if there is an
integer $r$ with $0\leqslant r\leqslant d$ and $0\leqslant p_0\leqslant\dots
\leqslant p_r\geqslant\dots\geqslant p_d\geqslant 0$.
An elementary but crucial property of reciprocal and unimodal polynomials
is the following.
\begin{lem}\label{lem:AqBq}
If $A(q)$ and $B(q)$ are reciprocal and unimodal polynomials, then so is
their product $A(q)B(q)$.
\end{lem}

Lemma \ref{lem:AqBq} is well known and its proof can be found, e.g.,
in \cite{Andrews75} or \cite[Proposition 1]{Stanley89}.


Similarly to the proof of \cite[Theorem 3.1]{GK}, the following lemma plays an important role in
the proof of Theorem \ref{thm:final}.
It is a slight generalization
of \cite[Proposition~10.1.(iii)]{ReSWAA},
which extracts the essentials out of an argument of Andrews \cite[Proof of Theorem~2]{AndrCB}.

\begin{lem} \label{lem:RSW}
Let $P(q)$ be a reciprocal and unimodal
polynomial and $m$ and $n$ positive integers with $m\leqslant n$.
Furthermore, assume that $A(q)=\frac {1-q^m} {1-q^n}P(q)$ is a polynomial
in $q$. Then $A(q)$ has non-negative coefficients.
\end{lem}
\pf See \cite[Lemma 7.1]{GK}.  \qed

\medskip
\noindent{\it Proof of Theorem {\rm\ref{thm:final}.}}
It is well known that the $q$-binomial coefficients are reciprocal and unimodal
polynomials in $q$ (cf.\ \cite[Ex.~7.75.d]{Stanley}), and by Lemma \ref{lem:AqBq}, so is
the product of two $q$-binomial coefficients. In view of Lemma~\ref{lem:RSW}, for proving Theorem~\ref{thm:final} it
is enough to show that the expression \eqref{eq:last} is a polynomial in $q$.
We shall accomplish this by a count of cyclotomic polynomials.

Recall the well-known fact that
$$
q^n-1=\prod _{d\mid n} ^{}\Phi_d(q),
$$
where $\Phi_d(q)$ denotes the $d$-th cyclotomic polynomial in $q$.
Consequently,
$$
\frac{1-q^{\gcd(am,m+n)}}{1-q^{m+n}}{am+bm-1\brack am}_q {an+bn\brack an}_q
=\prod _{d=2} ^{\min\{am+bm-1,\,an+bn\}}\Phi_d(q)^{e_d},
$$
with
\begin{align}
e_d&
=\chi(d\mid \gcd(am,m+n))-\chi(d\mid m+n)
+\fl{\frac {am+bm-1} {d}}
+\fl{\frac {an+bn} {d}} \notag\\
&\quad{}-\fl{\frac {am} {d}}
-\fl{\frac {bm-1} {d}}
-\fl{\frac {an} {d}}
-\fl{\frac {bn} {d}}, \label{eq:chi}
\end{align}
where $\chi(\mathcal S)=1$ if $\mathcal S$ is
true and $\chi(\mathcal S)=0$ otherwise.
This is clearly non-negative, unless
$d\mid m+n$ and $d\nmid \gcd(am,m+n)$.

So, let us assume that $d\mid m+n$ and $d\nmid \gcd(am,m+n)$, which means
that $d\nmid am$ and therefore
$$
\fl{\frac {am} {d}}=\fl{\frac {am-1} {d}}.
$$
Note that, when $d\mid m+n$, we have
\begin{align*}
\left\lfloor\frac{am+bm-1}{d}\right\rfloor+\left\lfloor\frac{an+bn}{d}\right\rfloor
&=\frac{(a+b)(m+n)}{d}-1, \\
\left\lfloor\frac{am-1}{d}\right\rfloor+\left\lfloor\frac{an}{d}\right\rfloor
&=\frac{a(m+n)}{d}-1,  \\
\left\lfloor\frac{bm-1}{d}\right\rfloor+\left\lfloor\frac{bn}{d}\right\rfloor
&=\frac{b(m+n)}{d}-1,
\end{align*}
and so $e_d=0$ is also non-negative in this case.
This completes the proof of polynomiality of \eqref{eq:last}.
\qed

\medskip
\noindent{\it Proof of Corollary {\rm\ref{cor:2}.}}
This follows immediately from Theorem~\ref{thm:final} and
the fact that $\gcd(am,m+n)\mid am$.  \qed

\section{Concluding remarks and open problems}
On January 2, 2014 T. Amdeberhan (personal communication) found the following generalization of
Theorem \ref{conj:sun3}. It would be interesting to give a proof of it.
\begin{conj}Let $a,b$ and $n$ be positive integers with $a>b$. Then
$$
(2bn+1)(2bn+3){2bn\choose bn}\left| 3(a-b)(3a-b){2an\choose an}{an\choose bn}.\right.
$$
\end{conj}

Let $[m]!=(1-q)\cdots (1-q^m)$. By a result of Warnaar and Zudilin \cite[Proposition 3]{WZ}, one sees that,
for any positive integer $n$, the polynomial
$$
\frac{[6n]![n]!}{[3n]![2n]!^2}
$$
has non-negative integer coefficients.
Similarly as before, we can prove the following generalizations of Theorem \ref{conj:sun3}
and the congruences \eqref{eq:number1} and \eqref{eq:number2}.
\begin{thm}\label{thm:remarks} Let $n$ be a positive integer. Then all of
\begin{align*}
&\frac{(1-q)[6n]![n]!}{(1-q^{2n+1})[3n]![2n]!^2}, \quad
\frac{(1-q^3)[6n]![n]!}{(1-q^{2n+3})[3n]![2n]!^2},\quad
\frac{(1-q)(1-q^3)[6n]![n]!}{(1-q^{2n+1})(1-q^{2n+3})[3n]![2n]!^2},\quad  \\[5pt]
&\frac{(1-q^3)(1-q^5)(1-q^7)[6n]![n]!}{(1-q^{2n+3})(1-q^{2n+5})(1-q^{2n+7})[3n]![2n]!^2}
\quad (n\geqslant 2),\\[5pt]
&\frac{(1-q^3)^2(1-q^5)(1-q^7)[6n]![n]!}{(1-q^{2n+1})(1-q^{2n+3})(1-q^{2n+5})(1-q^{2n+7})[3n]![2n]!^2}\quad
(n\geqslant 2)
\end{align*}
are polynomials in $q$.
\end{thm}

\begin{conj}
All the the polynomials in Theorem {\rm\ref{thm:remarks}} have non-negative integer coefficients.
\end{conj}

\begin{conj}\label{conj:unimodal}
Let $n\geqslant 2$. Then the polynomial $\frac{[6n]![n]!}{[3n]![2n]!^2}$ is unimodal.
\end{conj}

It is obvious that the polynomial $\frac{[6n]![n]!}{[3n]![2n]!^2}$ is reciprocal.
If Conjecture \ref{conj:unimodal} is true, then, applying Lemma \ref{lem:AqBq}, we conclude that the first two polynomials in
Theorem {\rm\ref{thm:remarks}} have non-negative integer coefficients.

It was conjectured by Warnaar and Zudilin (see \cite[Conjecture 1]{WZ}) that
$$
\frac{[15n]![2n]!}{[10n]![4n]![3n]!}
$$
has non-negative integer coefficients. Similarly, we have the following generalization
of Theorem \ref{conj:sun2}.
\begin{thm}\label{thm:remarks2} Let $n$ be a positive integer. Then both
\begin{align*}
&\frac{(1-q)[15n]![2n]!}{(1-q^{10n+1})[10n]![4n]![3n]!}, \quad\text{and}\quad
\frac{(1-q^3)(1-q^7)[15n]![2n]!}{(1-q)(1-q^{10n+3})[10n]![4n]![3n]!}  \label{eq:15n2n}
\end{align*}
are polynomials in $q$.
\end{thm}

\begin{conj}
The two polynomials in Theorem {\rm\ref{thm:remarks2}} have non-negative integer coefficients.
\end{conj}

\vskip 3mm
\noindent{\bf Acknowledgment.} The author thanks Quan-Hui Yang, T. Amdeberhan and V.H. Moll for helpful
comments on a previous version of this paper. This work was partially
supported by the Fundamental Research Funds for the Central Universities and
the National Natural Science Foundation of China (grant 11371144).

\renewcommand{\baselinestretch}{1}

\end{document}